\documentclass[preprint]{elsarticle}

\usepackage{lineno,hyperref}
\usepackage{mathtools,amsmath,amssymb}

\usepackage{graphicx}
\usepackage{tikz}

\journal{arXiv}

\usepackage{fullpage}

\newtheorem{theorem}{Theorem}%[section]
\newproof{pf}{Proof}

\makeatletter
\def\th@plain{%
  \thm@notefont{}% same as heading font
  \itshape % body font
}
\def\th@definition{%
  \thm@notefont{}% same as heading font
  \normalfont % body font
}
\makeatother
\DeclareFontFamily{U}{matha}{\hyphenchar\font45}
\DeclareFontShape{U}{matha}{m}{n}{ <5> <6> <7> <8> <9> <10> gen * matha <10.95> matha10 <12> <14.4> <17.28> <20.74> <24.88> matha12 }{}
\DeclareSymbolFont{matha}{U}{matha}{m}{n}
\DeclareFontFamily{U}{mathx}{\hyphenchar\font45}
\DeclareFontShape{U}{mathx}{m}{n}{ <5> <6> <7> <8> <9> <10> <10.95> <12> <14.4> <17.28> <20.74> <24.88> mathx10 }{}
\DeclareSymbolFont{mathx}{U}{mathx}{m}{n}
\DeclareMathSymbol{\obot} {2}{matha}{"6B}
\DeclareMathSymbol{\bigobot} {1}{mathx}{"CB}

\newcommand\mc{\mathcal}
\newcommand\mb{\mathbb}

\newcommand\ve{\varepsilon}%{\epsilon}
\newcommand\xt{\mathfrak{T}}

\DeclareMathOperator{\Span}{Span}

\DeclareMathOperator{\id}{Id}

\begin{document}

\begin{frontmatter}

\title{On the existence of a curvature tensor for given Jacobi operators}

%\author{Vladica Andreji\'c\fnref{myfootnote}}

\author[matf]{Vladica Andreji\'c}
%\author[matf]{V.~Andreji\'c}
\ead{andrew@matf.bg.ac.rs}
\author[matf]{Katarina Luki\'c}
\ead{katarina\_lukic@matf.bg.ac.rs}

\address[matf]{Faculty of Mathematics, University of Belgrade, Belgrade, Serbia}




\begin{abstract}
It is well known that the Jacobi operators completely determine the curvature tensor.
The question of existence of a curvature tensor for given Jacobi operators naturally arises, which is considered and solved in the previous work.
Unfortunately, although the published theorem is correct, its proof is incomplete because it contains some omissions, 
and the aim of this paper is to present a complete and accurate proof.
We also generalize the main theorem to the case of indefinite scalar product space.
Accordingly, we generalize the proportionality principle for Osserman algebraic curvature tensors.
\end{abstract}

\begin{keyword}
curvature tensor; Jacobi operators; duality principle; proportionality principle
\MSC[2020] Primary 53B20; Secondary 53B30
\end{keyword}

\end{frontmatter}

%\linenumbers

\section{Introduction}

The curvature is the most natural and most important invariant of Riemannian and pseudo-Riemannian geometry.
According to Osserman \cite{O}, the notion of curvature is one of the central concepts of differential geometry, 
distinguishing the geometric core of the subject from those aspects that are analytic, algebraic, or topological.
The curvature information is contained in the curvature tensor, which is difficult to work with, despite the many symmetries it possesses.
Extracting the geometrical information that is encoded therein is often quite a challenging task.
That is why Gromov \cite{Gro} described the curvature tensor as a little monster of (multi)linear algebra whose full geometric meaning remains obscure.
Therefore, instead of working with the curvature tensor itself, we often use Jacobi operators or sectional curvature
that are easier to handle and have a better geometrical interpretation.

Although the sectional curvature looks simpler than the curvature tensor, its importance arises from the fact that
knowledge of all sectional curvatures completely determines the curvature tensor (at points where the scalar product is known).
This refers to the Jacobi operators, which contain the very same information as the curvature tensor.
These well known results give the uniqueness of the curvature tensor and are purely algebraic in nature.
The question of existence of a curvature tensor for given Jacobi operators naturally arises, 
and as far as we know, the only work on the subject is \cite{A9}. 

Unfortunately, although the theorem published in \cite[Theorem 1]{A9} is correct, 
its proof is incomplete because it contains some omissions, and the aim of this paper is to present a complete and accurate proof.
This theorem of classical differential geometry is too important to allow it to remain incomplete or doubtful.
Moreover, the theorem is also valid in a pseudo-Riemannian case (the case of indefinite scalar product), and it is stated as follows.

\begin{theorem}\label{jacobiex}
Let $\mc{K}_X$ for all nonnull $X\in\mc{V}$ be a compatible family of self-adjoint endomorphisms on 
a (possibly indefinite) scalar product space $\mc{V}$ that satisfies $\mc{K}_XX=0$.
Then there exists a unique algebraic curvature tensor on $\mc{V}$ such that $\mc{K}_X$ are its Jacobi operators.
\end{theorem}

In Section \ref{PROOF1} we present the proof of Theorem \ref{jacobiex}.
In contrast to \cite{A9}, an algebraic curvature tensor is defined using formula \eqref{exist5} instead of \eqref{exist8}.
Although, the formulas \eqref{exist5}, \eqref{exist6}, and \eqref{exist8} are equivalent, 
the problem is that it is not clear why the derivative with respect to $t$ at $t=0$ in \eqref{exist6} exists, which is not justified in \cite{A9}.

On the one hand, although the condition \eqref{exist1} is written in the formulation of the theorem, in \cite{A9} it is not seen where it is used.
Here we give an elegant counterexample if the condition \eqref{exist1} does not hold, 
but what is much more important, we use \eqref{exist1} to prove the essential formula \eqref{exist4}.
On the other hand, the condition \eqref{exist2} in the formulation of the theorem in \cite{A9} is redundant, 
because it is a consequence of the other conditions.

However, the most important detail is that in \cite{A9} it was not checked whether the obtained algebraic curvature tensor really
has the Jacobi operators corresponding to the given family of endomorphisms.
This check is not trivial and requires proving the new crucial formula \eqref{exist4}. 

Finally, it is worth noting that Theorem \ref{jacobiex} is proved in the case of indefinite scalar product,
for which additional effort has been invested.
Accordingly, in Section \ref{APPLI}, the proportionality principle for Osserman algebraic curvature tensors is generalized to an indefinite case.

\section{The proof of Theorem \ref{jacobiex}}\label{PROOF1}

Let $(\mc{V},g)$ be a scalar product space and let $\ve_X=g(X,X)$ be the squared norm of $X\in\mc{V}$.
Since we want to cover pseudo-Riemannian settings, we allow the scalar product $g$ be indefinite, but it is certainly nondegenerate.
We say that a vector $X$ is null if $\ve_X=0$ and $X\neq 0$, while $X$ is nonnull if $\ve_X\neq 0$.

A tensor $R\in\xt^0_4(\mc{V})$ is said to be an {algebraic curvature tensor} on $(\mc{V},g)$ 
if it satisfies the usual $\mb{Z}_2$ symmetries as well as the first Bianchi identity.
More concretely, an algebraic curvature tensor $R$ has the following properties,
\begin{align}
&R(X,Y,Z,W)=-R(Y,X,Z,W),\label{tk1}\\
&R(X,Y,Z,W)=-R(X,Y,W,Z),\label{tk2}\\
&R(X,Y,Z,W)+R(Y,Z,X,W)+R(Z,X,Y,W)=0,\label{tk3}\\
&R(X,Y,Z,W)=R(Z,W,X,Y)\label{tk4},
\end{align}
that hold for all $X,Y,Z,W\in\mc{V}$, while it is well known that the equation \eqref{tk4} is a consequence of the first three equations. 
Raising the last index we obtain an algebraic curvature operator $\mc{R}=R^\sharp\in\xt^1_3(\mc{V})$.

The Jacobi operator is a linear operator $\mc{J}_X\colon\mc{V}\to\mc{V}$ 
defined for each $X\in\mc{V}$ by $\mc{J}_X(Y)=\mc{R}(Y,X)X$ for all $Y\in\mc{V}$.
The symmetries of $R$ show that each Jacobi operator is self-adjoint, and any two are mutually compatible, which means that 
\begin{equation}\label{jak2}
g(\mc{J}_XY,Y)=g(\mc{J}_YX,X)
\end{equation}
holds for all $X,Y\in\mc{V}$.
Thus, the Jacobi operators are self-adjoint endomorphisms on $\mc{V}$, which satisfy the compatibility condition \eqref{jak2}.
Using the symmetries of $R$, the straightforward calculations for arbitrary $X,Y,Z,W\in\mc{V}$ give
\begin{equation*}
\begin{aligned}
3R(X,Y,Z,W)&=R(X,Y,Z,W)-R(X,Y,W,Z)+(-R(Y,Z,X,W)-R(Z,X,Y,W))\\
&=(R(X,Y,Z,W)+R(X,Z,Y,W))-(R(X,Y,W,Z)+R(X,W,Y,Z))\\
&=g((\mc{J}_{Y+Z}-\mc{J}_{Y}-\mc{J}_{Z})X,W)-g((\mc{J}_{Y+W}-\mc{J}_{Y}-\mc{J}_{W})X,Z)%\\
\end{aligned}
\end{equation*} 
and therefore
\begin{equation}\label{exist0}
3R(X,Y,Z,W)=g((\mc{J}_{Y+Z}-\mc{J}_{Y}-\mc{J}_{Z})W-(\mc{J}_{Y+W}-\mc{J}_{Y}-\mc{J}_{W})Z,X).
\end{equation} 
The formula \eqref{exist0} shows that the Jacobi operators completely determine the curvature tensor
(moreover, the sectional curvatures do this as well), which is well known.

The question of existence of a curvature tensor for given Jacobi operators naturally arises. 
Let us suppose that we know self-adjoint endomorphisms $\mc{K}_X$ on $\mc{V}$ for each nonnull $X\in\mc{V}$,
such that they are compatible in the sense that \eqref{jak2} holds.
Is there an algebraic curvature tensor $R$ on $(\mc{V},g)$, such that $\mc{J}_X=\mc{K}_X$ holds for all nonnull $X\in\mc{V}$?

Consider the condition $\mc{J}_XX=0$ which holds for all Jacobi operators. 
However, the construction $\mc{K}_X=\ve_X\id$ gives self-adjoint endomorphisms on $\mc{V}$ that are compatible,
$g(\ve_X\id Y,Y)=\ve_X\ve_Y=g(\ve_Y\id X,X)$, but $\ve_X\id X=\ve_XX\neq 0$ holds for a nonnull $X\in\mc{V}$.
Therefore, we add the natural condition that 
\begin{equation}\label{exist1}
\mc{K}_XX=0,
\end{equation}
holds for any nonnull $X\in\mc{V}$. 

The first step is to extend the family $\mc{K}_X$ for all $X\in\mc{V}$.
The natural extension is $\mc{K}_0=0$, which completes the family in the Riemannian setting.
However, if the scalar product is indefinite, then we need to define $\mc{K}_X$ for any null $X\in\mc{V}$.

Let $X,Y,X+Y,X-Y\in\mc{V}$ be nonnull. Using \eqref{jak2}, from 
\begin{equation*}
\begin{aligned}
g(\mc{K}_{X\pm Y}Z,Z)=g(\mc{K}_Z(X\pm Y),X\pm Y)
&=g(\mc{K}_Z X,X) \pm 2g(\mc{K}_Z X, Y) +g(\mc{K}_Z Y,Y)\\
&=g(\mc{K}_X Z,Z) \pm 2g(\mc{K}_Z X, Y) +g(\mc{K}_Y Z,Z),
\end{aligned}
\end{equation*}
it follows $g(\mc{K}_{X+Y}Z,Z)+g(\mc{K}_{X-Y}Z,Z)=2g(\mc{K}_X Z,Z)+2g(\mc{K}_Y Z,Z)$ for any nonnull $Z\in\mc{V}$.
The polarization $Z=V+W$ gives $g((\mc{K}_{X+Y}+\mc{K}_{X-Y}-2\mc{K}_X-2\mc{K}_Y)V,W)=0$, whenever $\ve_V\ve_W\ve_{V+W}\neq 0$. 
If $(E_1,E_2,\dots,E_n)$ is an orthonormal basis in $(\mc{V},g)$, then we can consider the orthogonal basis $(E_1,2E_2,\dots,nE_n)$,
which provides $\ve_{iE_i}=i^2\ve_{E_i}\neq 0$ and $\ve_{iE_i+jE_j}=i^2\ve_{E_i}+j^2\ve_{E_j}\neq 0$ for $1\leq i\neq j\leq n$.
Hence,
\begin{equation}\label{existnull}
\mc{K}_{X+Y}+\mc{K}_{X-Y}=2\mc{K}_X+2\mc{K}_Y
\end{equation}
holds, whenever $\ve_X\ve_Y\ve_{X+Y}\ve_{X-Y}\neq 0$.

The equation \eqref{existnull} motivates us to define $\mc{K}_N$ for a null $N\in\mc{V}$ by
\begin{equation*}
2\mc{K}_N=\mc{K}_{N+X}+\mc{K}_{N-X}-2\mc{K}_X,
\end{equation*}
whenever the right hand side is defined, which immediately shows that $\mc{K}_N$ is a self-adjoint endomorphism on $\mc{V}$. 
Because of \eqref{existnull}, if $\ve_X\ve_Y\ve_{N+X}\ve_{N-X}\ve_{N+Y}\ve_{N-Y}\neq 0$, then we have
\begin{equation*}
\begin{aligned}
2(\mc{K}_{N+X}+\mc{K}_{N-X}-2\mc{K}_X)&=(\mc{K}_{N+X+Y}+\mc{K}_{N+X-Y}-2\mc{K}_Y)+(\mc{K}_{N-X+Y}+\mc{K}_{N-X-Y}-2\mc{K}_Y)-4\mc{K}_X\\
&=(\mc{K}_{N+X+Y}+\mc{K}_{N-X+Y}-2\mc{K}_X)+(\mc{K}_{N+X-Y}+\mc{K}_{N-X-Y}-2\mc{K}_X)-4\mc{K}_Y\\
&=2(\mc{K}_{N+Y}+\mc{K}_{N-Y}-2\mc{K}_Y),
\end{aligned}
\end{equation*}
whenever $\ve_{N+X+Y}\ve_{N+X-Y}\ve_{N-X+Y}\ve_{N-X-Y}\neq 0$. Otherwise, we can use 
\begin{equation*}
\mc{K}_{N+X}+\mc{K}_{N-X}-2\mc{K}_X= \mc{K}_{N+tX}+\mc{K}_{N-tX}-2\mc{K}_{tX} =\mc{K}_{N+Y}+\mc{K}_{N-Y}-2\mc{K}_Y,
\end{equation*}
where $t$ is not a root of 
\begin{equation*}
\ve_{tX}\ve_{N+tX}\ve_{N-tX}\ve_{N+X+tX}\ve_{N-X+tX}\ve_{N+X-tX}\ve_{N-X-tX}\ve_{N+Y+tX}\ve_{N-Y+tX}\ve_{N+Y-tX}\ve_{N-Y-tX}=0,
\end{equation*}
which is a polynomial equation of degree $22$. 
This proves that $\mc{K}_N$ does not depend on the choice of $X$, and therefore $\mc{K}_N$ is well-defined.

If we use $Z\in\mc{V}$ that satisfies $\ve_{N+Z}\ve_{N-Z}\ve_{Z}\neq 0$, then the equation 
\begin{equation*}
\begin{aligned}
2g(\mc{K}_NX,X)&=g((\mc{K}_{N+Z}+\mc{K}_{N-Z}-2\mc{K}_Z)X,X)\\
&=g(\mc{K}_X(N+Z),N+Z)+g(\mc{K}_X(N-Z),N-Z)-2g(\mc{K}_XZ,Z)=2g(\mc{K}_XN,N)
\end{aligned}
\end{equation*}
holds for any null $N$ and any nonnull $X$, which means that endomorphisms $\mc{K}_N$ and $\mc{K}_X$ are compatible.
With that in mind, we see that the very same equation holds when both $N$ and $X$ are null.
In this way, we obtain $\mc{K}_X$ for any $X\in\mc{V}$, and this extended family remains compatible. 

Let $N\in\mc{V}$ be null. For any $X\in\mc{V}$ such that $\ve_X\ve_{N+X}\ve_{N-X}\neq 0$ holds, 
using the properties \eqref{exist1} and \eqref{jak2} we have
\begin{equation*}
\begin{aligned}
2g(\mc{K}_NN,X)&=g(\mc{K}_{N+X}N+\mc{K}_{N-X}N-2\mc{K}_X N,X)=-g(\mc{K}_{N+X}X,X)+g(\mc{K}_{N-X}X,X)-2g(\mc{K}_X X,N)\\
&=-g(\mc{K}_{X}(N+X),N+X)+g(\mc{K}_{X}(N-X),N-X)=0.
\end{aligned}
\end{equation*}
If $(E_1,E_2,\dots,E_n)$ is an orthonormal basis in $(\mc{V},g)$, then we create an orthogonal basis $(mE_1,mE_2,\dots,mE_n)$,
where $m>\max_{1\leq i\leq n} {|2g(N,E_i)|}$ to provide $\ve_{mE_i}\neq 0$ and $\ve_{N\pm mE_i}\neq 0$ for $1\leq i\leq n$,
which yields $K_NN=0$.

The second step considers a compatible family of self-adjoint endomorphisms $\mc{K}_X$ for all $X\in\mc{V}$ such that $\mc{K}_XX=0$ holds.
For all $X,Y\in\mc{V}$ and $t\in\mb{R}$ we have
%\begin{equation*}
$g(\mc{K}_{tX}Y,Y)=g(\mc{K}_{Y}tX,tX)=t^2g(\mc{K}_{Y}X,X)=g(t^2\mc{K}_{X}Y,Y)$,
%\end{equation*}
where the polarization $Y=V+W$ gives $g(\mc{K}_{tX}V,W)=g(t^2\mc{K}_{X}V,W)$ for $V,W\in\mc{V}$, 
so since $g$ is nondegenerate, it follows
\begin{equation}\label{exist2}
\mc{K}_{tX}=t^2\mc{K}_{X},
\end{equation}
which is a natural property of Jacobi operators.

For all $X,Y,Z\in\mc{V}$ and $t\in\mb{R}$, using the compatibility we have
\begin{equation*}
\begin{aligned}
g(\mc{K}_{X+tY}Z,Z)&=g(\mc{K}_Z(X+tY),X+tY)=g(\mc{K}_ZX,X)+2tg(\mc{K}_ZX,Y)+t^2g(\mc{K}_ZY,Y)\\
&=g(\mc{K}_XZ,Z)+2tg(\mc{K}_ZX,Y)+t^2g(\mc{K}_YZ,Z),
\end{aligned}
\end{equation*}
which implies
\begin{equation*}
g((\mc{K}_{X+tY}-\mc{K}_X-t^2\mc{K}_Y )Z,Z)=2tg(\mc{K}_ZX,Y)=tg((\mc{K}_{X+Y}-\mc{K}_X-\mc{K}_Y )Z,Z).
\end{equation*}
After the polarization $Z=V+W$ we get
%\begin{equation*}
$g((\mc{K}_{X+tY}-\mc{K}_X-t^2\mc{K}_Y )V,W)=tg((\mc{K}_{X+Y}-\mc{K}_X-\mc{K}_Y )V,W)$ for all $V,W\in\mc{V}$,
%\end{equation*}
and therefore, since $g$ is nondegenerate, we obtain a generalization of \eqref{existnull},
\begin{equation}\label{exist3}
\mc{K}_{X+tY}-\mc{K}_X-t^2\mc{K}_Y =t( \mc{K}_{X+Y}-\mc{K}_X-\mc{K}_Y ),
\end{equation}
which yields
%\begin{equation*}
$\mc{K}_{X+tY}= t\mc{K}_{X+Y} +(1-t)\mc{K}_X +(t^2-t)\mc{K}_Y$.
%\end{equation*}
Using the property \eqref{exist1}, we get
\begin{equation*}
\begin{aligned}
0=\mc{K}_{X+tY}(X+tY)&= t\mc{K}_{X+Y}(X+tY) +(1-t)\mc{K}_X(X+tY) +(t^2-t)\mc{K}_Y(X+tY)\\
&= t\mc{K}_{X+Y}((t-1)Y) +(1-t)\mc{K}_X(tY) + t(t-1)\mc{K}_Y(X)\\
&= t(t-1)(\mc{K}_{X+Y}Y -\mc{K}_XY +\mc{K}_YX),
\end{aligned}
\end{equation*}
which for $t\in\mb{R}\setminus\{0,1\}$ implies
\begin{equation}\label{exist4}
\mc{K}_{X+Y}Y -\mc{K}_XY +\mc{K}_YX=0
\end{equation}
for all $X,Y\in\mc{V}$.

Let us use a compatible family of self-adjoint endomorphisms $\mc{K}_X$ on $\mc{V}$ that satisfies $\mc{K}_XX=0$ 
to define a map $R\colon \mc{V}^4\to\mb{R}$, resembling the formula \eqref{exist0} by
\begin{equation}\label{exist5}
3R(X,Y,Z,W)=g((\mc{K}_{Y+Z}-\mc{K}_{Y}-\mc{K}_{Z})W-(\mc{K}_{Y+W}-\mc{K}_{Y}-\mc{K}_{W})Z,X),
\end{equation} 
for all $X,Y,Z,W\in\mc{V}$.
From \eqref{exist3}, if we take the limit where $t$ tends to zero, then it follows
\begin{equation*}
g((\mc{K}_{Y+Z}-\mc{K}_{Y}-\mc{K}_{Z})W,X)=\frac{g(\mc{K}_{Y+tZ}W,X)-g(\mc{K}_{Y}W,X)}{t}-tg(\mc{K}_{Z}W,X)
=\left.\frac{\partial}{\partial t}\right|_{t=0}  g(\mc{K}_{Y+tZ}W,X),
\end{equation*}
and therefore the equation \eqref{exist5} is equivalent to 
\begin{equation}\label{exist6}
3R(X,Y,Z,W)= \left.\frac{\partial}{\partial t}\right|_{t=0}  g(\mc{K}_{Y+tZ}W-\mc{K}_{Y+tW}Z ,X).
\end{equation} 
However, from 
\begin{equation*}
\begin{aligned}
2g(\mc{K}_{Y+tZ}W,X)&=\left.\frac{\partial}{\partial s}\right|_{s=0}\Big( g(\mc{K}_{Y+tZ}X,X) +2sg(\mc{K}_{Y+tZ}W,X) + s^2g(\mc{K}_{Y+tZ}W,W) \Big)\\
&=\left.\frac{\partial}{\partial s}\right|_{s=0}\Big( g(\mc{K}_{Y+tZ}W,W) +2sg(\mc{K}_{Y+tZ}W,X) + s^2g(\mc{K}_{Y+tZ}X,X) \Big),
\end{aligned}
\end{equation*} 
we obtain
\begin{equation}\label{exist7}
2g(\mc{K}_{Y+tZ}W,X)=\left.\frac{\partial}{\partial s}\right|_{s=0}\mu(X+sW,Y+tZ)=\left.\frac{\partial}{\partial s}\right|_{s=0}\mu(W+sX,Y+tZ),
\end{equation} 
where $\mu(X,Y)=g(\mc{K}_YX,X)=g(\mc{K}_XY,Y)$ for all $X,Y\in\mc{V}$.
The equality on the left side in \eqref{exist7} shows that \eqref{exist6} is equivalent to
\begin{equation}\label{exist8}
6R(X,Y,Z,W)= \left.\frac{\partial^2}{\partial s \partial t}\right|_{s=0,t=0}(\mu(X+sW,Y+tZ) -\mu(X+sZ,Y+tW)).
\end{equation} 
In this way, we obtain the formula \eqref{exist8}, which is known as another proof that the sectional curvature
completely determines the curvature tensor, see Lee \cite[Proposition 13.27]{Lej}.

The definition \eqref{exist5} of $R$ is equivalent to \eqref{exist8}, but the latter is easier to prove that $R$
is an algebraic curvature tensor. 
The property \eqref{tk2} follows directly from \eqref{exist8}, while \eqref{tk1} is a consequence of 
commutativity $\partial_s\partial_t=\partial_t\partial_s$. The equality on the right side in \eqref{exist7}
helps us to easily see \eqref{tk3}.
From \eqref{exist6}, $R(X,Y,Z,W)$ is obviously linear by $X$, but due to the already proven symmetries, 
where \eqref{tk4} automatically follows, $R$ is multi-linear, which proves that $R$ is an algebraic curvature tensor.

It remains to show that $\mc{K}_X$ for $X\in\mc{V}$ are the Jacobi operators for $R$.
From the definition \eqref{exist5} of $R$, using $g(\mc{J}_YW,X)=R(W,Y,Y,X)$ we have
%\begin{equation*}
$3g(\mc{J}_YW,X)=g(2\mc{K}_{Y}W - \mc{K}_{Y+W}Y +\mc{K}_{W}Y,X)$, which implies
%\end{equation*} 
\begin{equation*}
3\mc{J}_YW=2\mc{K}_{Y}W - \mc{K}_{Y+W}Y +\mc{K}_{W}Y.
\end{equation*} 
According to \eqref{exist4} we have $-\mc{K}_{Y+W}Y=\mc{K}_{Y+W}W=\mc{K}_YW -\mc{K}_WY$,
and therefore $3\mc{J}_YW=3\mc{K}_{Y}W$ for all $Y,W\in\mc{V}$, which gives $\mc{J}_Y=\mc{K}_Y$.
This finally proves Theorem \ref{jacobiex}.

\section{Applications to Osserman manifolds}\label{APPLI}

Theorem \ref{jacobiex} was basically invented in \cite{A9} for application to Osserman curvature tensors.
The aim of this section is to slightly improve things from \cite{A9} and give a generalization in pseudo-Riemannian settings.

Let $R$ be an algebraic curvature tensor on a (possibly indefinite) scalar product space $(\mc{V},g)$.
Since $\mc{J}_XX=0$ and $g(\mc{J}_XY,X)=0$, the Jacobi operator $\mc{J}_X$ for a nonnull $X\in\mc{V}$ 
is completely determined by its restriction $\widetilde{\mc{J}}_X\colon X^\perp \to X^\perp$ called the reduced Jacobi operator.

We say that $R$ is Osserman if the polynomial $\det(\lambda\id-\mc{J}_X/\ve_X)$ is independent of a nonnull $X\in\mc{V}$.
More generally, we say that $R$ is $k$-root if $\widetilde{\mc{J}}_X$ has exactly $k$ distinct eigenvalues
(counting complex roots) for any nonnull $X\in\mc{V}$.
If for each nonnull $X\in\mc{V}$ there exists an orthonormal eigenbasis in $\mc{V}$ related to $\mc{J}_X$, we say that $R$ is Jacobi-diagonalizable.
In the case of Jacobi-diagonalizable $k$-root Osserman $R$, any nonnull $X\in\mc{V}$ 
allows the spectral decomposition $\mc{V}=\bigoplus_{i=0}^k \mc{V}_i(X)$ 
as the orthogonal direct sum of the (generalized) eigenspaces $\mc{V}_i(X)=\ker(\widetilde{\mc{J}}_X-\ve_X\lambda_i\id)$ for $1\leq i\leq k$, 
with additional $\mc{V}_0(X)=\Span\{X\}$.
We say that a Jacobi-diagonalizable $k$-root Osserman $R$ is Jacobi-proportional if for any pair of nonnull vectors $X,Y\in\mc{V}$ we have the proportionality
\begin{equation*}
\ve_X(\ve_{Y_0}, \ve_{Y_1}, \dots, \ve_{Y_k})=\ve_Y(\ve_{X_0}, \ve_{X_1}, \dots, \ve_{X_k}),
\end{equation*}
where $X=\sum_{i=0}^kX_i$ and $Y=\sum_{j=0}^kY_j$ are decomposed such that $X_i\in\mc{V}_i(Y)$ and $Y_j\in\mc{V}_j(X)$.

%%%%Picture
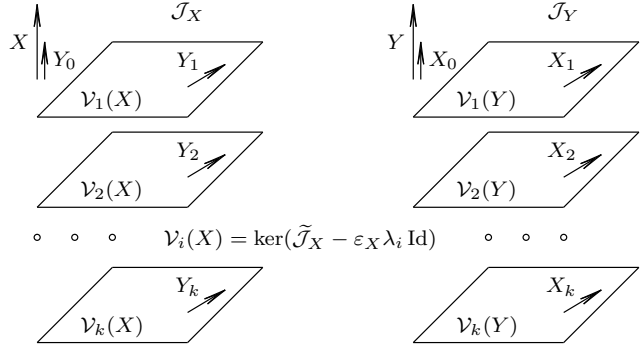
\begin{center}
%\input{X101.tkz}
%\end{center}
\begin{tikzpicture}
\clip (0,0) rectangle (10.000000,4.700000);
{\footnotesize

% Drawing segment D1 D2
\draw [line width=0.016cm] (1.000000,0.100000) -- (2.000000,1.100000);%

% Drawing segment D3 D2
\draw [line width=0.016cm] (4.000000,1.100000) -- (2.000000,1.100000);%

% Drawing segment D3 D4
\draw [line width=0.016cm] (4.000000,1.100000) -- (3.000000,0.100000);%

% Drawing segment D1 D4
\draw [line width=0.016cm] (1.000000,0.100000) -- (3.000000,0.100000);%

% Printing text at point Q
\draw (1.470000,0.070000) node [anchor=south west] { $\mc{V}_k(X)$ };%

% Drawing segment C1 C2
\draw [line width=0.016cm] (1.000000,1.900000) -- (2.000000,2.900000);%

% Drawing segment C3 C2
\draw [line width=0.016cm] (4.000000,2.900000) -- (2.000000,2.900000);%

% Drawing segment C3 C4
\draw [line width=0.016cm] (4.000000,2.900000) -- (3.000000,1.900000);%

% Drawing segment C1 C4
\draw [line width=0.016cm] (1.000000,1.900000) -- (3.000000,1.900000);%

% Printing text at point Q
\draw (1.470000,1.870000) node [anchor=south west] { $\mc{V}_2(X)$ };%

% Marking point P by circle
\draw [line width=0.016cm] (1.000000,1.500000) circle (0.040000);%

% Marking point P by circle
\draw [line width=0.016cm] (1.500000,1.500000) circle (0.040000);%

% Marking point P by circle
\draw [line width=0.016cm] (2.000000,1.500000) circle (0.040000);%

% Printing text at point P
\draw (4.500000,1.500000) node  { $\mc{V}_i(X)=\ker(\widetilde{\mc{J}}_X-\ve_X\lambda_i\id)$ };%

% Drawing segment B1 B2
\draw [line width=0.016cm] (1.000000,3.100000) -- (2.000000,4.100000);%

% Drawing segment B3 B2
\draw [line width=0.016cm] (4.000000,4.100000) -- (2.000000,4.100000);%

% Drawing segment B3 B4
\draw [line width=0.016cm] (4.000000,4.100000) -- (3.000000,3.100000);%

% Drawing segment B1 B4
\draw [line width=0.016cm] (1.000000,3.100000) -- (3.000000,3.100000);%

% Printing text at point Q
\draw (1.470000,3.070000) node [anchor=south west] { $\mc{V}_1(X)$ };%

% Drawing vector A1 A2
\draw [line width=0.016cm] (1.000000,3.600000) -- (1.000000,4.600000);%
\draw [line width=0.016cm] (0.960842,4.302567) -- (1.000000,4.600000);%
\draw [line width=0.016cm] (0.960842,4.302567) -- (1.000000,4.500000);%
\draw [line width=0.016cm] (1.039158,4.302567) -- (1.000000,4.600000);%
\draw [line width=0.016cm] (1.039158,4.302567) -- (1.000000,4.500000);%

% Printing text at point M
\draw (1.000000,4.100000) node [anchor=east] { $X$ };%

% Printing text at point T
\draw (3.000000,4.500000) node  { $\mc{J}_X$ };%

% Drawing vector H1 H2
\draw [line width=0.016cm] (6.000000,3.600000) -- (6.000000,4.600000);%
\draw [line width=0.016cm] (5.960842,4.302567) -- (6.000000,4.600000);%
\draw [line width=0.016cm] (5.960842,4.302567) -- (6.000000,4.500000);%
\draw [line width=0.016cm] (6.039158,4.302567) -- (6.000000,4.600000);%
\draw [line width=0.016cm] (6.039158,4.302567) -- (6.000000,4.500000);%

% Printing text at point M
\draw (6.000000,4.100000) node [anchor=east] { $Y$ };%

% Printing text at point T
\draw (8.000000,4.500000) node  { $\mc{J}_Y$ };%

% Drawing segment E1 E2
\draw [line width=0.016cm] (6.000000,3.100000) -- (7.000000,4.100000);%

% Drawing segment E3 E2
\draw [line width=0.016cm] (9.000000,4.100000) -- (7.000000,4.100000);%

% Drawing segment E3 E4
\draw [line width=0.016cm] (9.000000,4.100000) -- (8.000000,3.100000);%

% Drawing segment E1 E4
\draw [line width=0.016cm] (6.000000,3.100000) -- (8.000000,3.100000);%

% Printing text at point Q
\draw (6.470000,3.070000) node [anchor=south west] { $\mc{V}_1(Y)$ };%

% Drawing segment F1 F2
\draw [line width=0.016cm] (6.000000,1.900000) -- (7.000000,2.900000);%

% Drawing segment F3 F2
\draw [line width=0.016cm] (9.000000,2.900000) -- (7.000000,2.900000);%

% Drawing segment F3 F4
\draw [line width=0.016cm] (9.000000,2.900000) -- (8.000000,1.900000);%

% Drawing segment F1 F4
\draw [line width=0.016cm] (6.000000,1.900000) -- (8.000000,1.900000);%

% Printing text at point Q
\draw (6.470000,1.870000) node [anchor=south west] { $\mc{V}_2(Y)$ };%

% Drawing segment G1 G2
\draw [line width=0.016cm] (6.000000,0.100000) -- (7.000000,1.100000);%

% Drawing segment G3 G2
\draw [line width=0.016cm] (9.000000,1.100000) -- (7.000000,1.100000);%

% Drawing segment G3 G4
\draw [line width=0.016cm] (9.000000,1.100000) -- (8.000000,0.100000);%

% Drawing segment G1 G4
\draw [line width=0.016cm] (6.000000,0.100000) -- (8.000000,0.100000);%

% Printing text at point Q
\draw (6.470000,0.070000) node [anchor=south west] { $\mc{V}_k(Y)$ };%

% Marking point P by circle
\draw [line width=0.016cm] (7.000000,1.500000) circle (0.040000);%

% Marking point P by circle
\draw [line width=0.016cm] (7.500000,1.500000) circle (0.040000);%

% Marking point P by circle
\draw [line width=0.016cm] (8.000000,1.500000) circle (0.040000);%

% Drawing vector Z1 Z2
\draw [line width=0.016cm] (1.100000,3.600000) -- (1.100000,4.100000);%
\draw [line width=0.016cm] (1.060842,3.802567) -- (1.100000,4.100000);%
\draw [line width=0.016cm] (1.060842,3.802567) -- (1.100000,4.000000);%
\draw [line width=0.016cm] (1.139158,3.802567) -- (1.100000,4.100000);%
\draw [line width=0.016cm] (1.139158,3.802567) -- (1.100000,4.000000);%

% Printing text at point M
\draw (1.100000,3.850000) node [anchor=west] { $Y_0$ };%

% Drawing vector Z1 Z2
\draw [line width=0.016cm] (3.000000,0.500000) -- (3.500000,0.800000);%
\draw [line width=0.016cm] (3.224806,0.680549) -- (3.500000,0.800000);%
\draw [line width=0.016cm] (3.224806,0.680549) -- (3.414251,0.748550);%
\draw [line width=0.016cm] (3.265099,0.613394) -- (3.500000,0.800000);%
\draw [line width=0.016cm] (3.265099,0.613394) -- (3.414251,0.748550);%

% Printing text at point M
\draw (3.280000,0.620000) node [anchor=south east] { $Y_k$ };%

% Drawing vector Z1 Z2
\draw [line width=0.016cm] (3.000000,2.300000) -- (3.500000,2.600000);%
\draw [line width=0.016cm] (3.224806,2.480549) -- (3.500000,2.600000);%
\draw [line width=0.016cm] (3.224806,2.480549) -- (3.414251,2.548550);%
\draw [line width=0.016cm] (3.265099,2.413394) -- (3.500000,2.600000);%
\draw [line width=0.016cm] (3.265099,2.413394) -- (3.414251,2.548550);%

% Printing text at point M
\draw (3.280000,2.420000) node [anchor=south east] { $Y_2$ };%

% Drawing vector Z1 Z2
\draw [line width=0.016cm] (3.000000,3.500000) -- (3.500000,3.800000);%
\draw [line width=0.016cm] (3.224806,3.680549) -- (3.500000,3.800000);%
\draw [line width=0.016cm] (3.224806,3.680549) -- (3.414251,3.748550);%
\draw [line width=0.016cm] (3.265099,3.613394) -- (3.500000,3.800000);%
\draw [line width=0.016cm] (3.265099,3.613394) -- (3.414251,3.748550);%

% Printing text at point M
\draw (3.280000,3.620000) node [anchor=south east] { $Y_1$ };%

% Drawing vector Z1 Z2
\draw [line width=0.016cm] (6.100000,3.600000) -- (6.100000,4.100000);%
\draw [line width=0.016cm] (6.060842,3.802567) -- (6.100000,4.100000);%
\draw [line width=0.016cm] (6.060842,3.802567) -- (6.100000,4.000000);%
\draw [line width=0.016cm] (6.139158,3.802567) -- (6.100000,4.100000);%
\draw [line width=0.016cm] (6.139158,3.802567) -- (6.100000,4.000000);%

% Printing text at point M
\draw (6.100000,3.850000) node [anchor=west] { $X_0$ };%

% Drawing vector Z1 Z2
\draw [line width=0.016cm] (8.000000,0.500000) -- (8.500000,0.800000);%
\draw [line width=0.016cm] (8.224806,0.680549) -- (8.500000,0.800000);%
\draw [line width=0.016cm] (8.224806,0.680549) -- (8.414251,0.748550);%
\draw [line width=0.016cm] (8.265099,0.613394) -- (8.500000,0.800000);%
\draw [line width=0.016cm] (8.265099,0.613394) -- (8.414251,0.748550);%

% Printing text at point M
\draw (8.280000,0.620000) node [anchor=south east] { $X_k$ };%

% Drawing vector Z1 Z2
\draw [line width=0.016cm] (8.000000,2.300000) -- (8.500000,2.600000);%
\draw [line width=0.016cm] (8.224806,2.480549) -- (8.500000,2.600000);%
\draw [line width=0.016cm] (8.224806,2.480549) -- (8.414251,2.548550);%
\draw [line width=0.016cm] (8.265099,2.413394) -- (8.500000,2.600000);%
\draw [line width=0.016cm] (8.265099,2.413394) -- (8.414251,2.548550);%

% Printing text at point M
\draw (8.280000,2.420000) node [anchor=south east] { $X_2$ };%

% Drawing vector Z1 Z2
\draw [line width=0.016cm] (8.000000,3.500000) -- (8.500000,3.800000);%
\draw [line width=0.016cm] (8.224806,3.680549) -- (8.500000,3.800000);%
\draw [line width=0.016cm] (8.224806,3.680549) -- (8.414251,3.748550);%
\draw [line width=0.016cm] (8.265099,3.613394) -- (8.500000,3.800000);%
\draw [line width=0.016cm] (8.265099,3.613394) -- (8.414251,3.748550);%

% Printing text at point M
\draw (8.280000,3.620000) node [anchor=south east] { $X_1$ };%
}
\end{tikzpicture}

\end{center}

The main idea from \cite{A9} was to consider a Jacobi-proportional Jacobi-diagonalizable Osserman algebraic curvature tensor $R$,
and use $\mc{J}_XV=\ve_X\lambda_iV$ to define $\mc{K}_XV=\ve_X\mu_iV$ for arbitrary $\mu_1,\dots,\mu_k\in\mb{R}$.
Of course, we add $\mc{K}_XX=0$, which means $\mu_0=0$.
In this way, we obtain self-adjoint endomorphisms $\mc{K}_X$ for any nonnull $X\in\mc{V}$ 
by keeping the existing eigenspaces and replacing the eigenvalues. 

Since $R$ is Jacobi-proportional, if we decompose nonnull $X,Y\in\mc{V}$ by $X=\sum_{i=0}^kX_i$ and $Y=\sum_{j=0}^kY_j$ 
such that $X_i\in\mc{V}_i(Y)$ and $Y_j\in\mc{V}_j(X)$ hold, then we have $\ve_X\ve_{Y_i}=\ve_Y\ve_{X_i}$ for all $0\leq i\leq k$.
Hence, we have
\begin{equation*}
\begin{aligned}
g(\mc{K}_X Y,Y)&=g\left(\sum_{i=0}^k\ve_X\mu_iY_i, Y\right)=\ve_X\sum_{i=0}^k\mu_i g({Y_i},Y)=\ve_X\sum_{i=0}^k\mu_i\ve_{Y_i},\\
g(\mc{K}_Y X,X)&=g\left(\sum_{i=0}^k\ve_Y\mu_iX_i, X\right)=\ve_Y\sum_{i=0}^k\mu_i g({X_i},X)=\ve_Y\sum_{i=0}^k\mu_i\ve_{X_i},
\end{aligned}
\end{equation*}
and therefore $g(\mc{K}_Y X,X)=g(\mc{K}_X Y,Y)$ holds for nonnull $X,Y\in\mc{V}$.
Thus, we fulfilled the conditions of Theorem \ref{jacobiex}, which gives a generalization of the theorem from \cite[Theorem 2]{A9}.

\begin{theorem}\label{newoss}
If there exists a Jacobi-proportional Jacobi-diagonalizable $k$-root Osserman algebraic curvature tensor such that 
$\ker(\lambda\id-\widetilde{\mc{J}}_X)=\prod_{i=1}^k(\lambda-\ve_X\lambda_i)^{\nu_i}$,
then for any scalars $\mu_1,\dots,\mu_k\in\mb{R}$ there is a new Osserman algebraic curvature tensor such that
$\ker(\lambda\id-\widetilde{\mc{J}}_X)=\prod_{i=1}^k(\lambda-\ve_X\mu_i)^{\nu_i}$.
\end{theorem}

We say that an algebraic curvature tensor $R$ is semi-Clifford if $R=\mu_0R^1+\sum_{i=1}^m\mu_iR^{J_i}$
for some anti-commutative family of skew-adjoint complex or product structures $J_i$ for $1\leq i\leq m$,
and $\mu_0,\dots,\mu_m\in\mb{R}$, where
\begin{equation*}
\begin{aligned}
R^1(X,Y,Z,W)&=g(Y,Z)g(X,W)-g(X,Z)g(Y,W),\\
R^{J_i}(X,Y,Z,W)&= g(J_iX,Z)g(J_iY,W) -g(J_iY,Z)g(J_iX,W) +2g(J_iX,Y)g(J_iZ,W),
\end{aligned}
\end{equation*}
holds for $X,Y,Z,W\in\mc{V}$.
It is well known that any semi-Clifford $R$ is Jacobi-diagonalizable Osserman.
If all $J_i$ are complex structures ($J_i^2=-\id$) then we say that $R$ is Clifford. 
Let us remark that product structures ($J_i^2=\id$) exist only in the case of neutral signature,
so if the signature is not neutral, any semi-Clifford $R$ is Clifford.

The main application of Theorem \ref{newoss} concerns the Osserman conjecture in the Riemannian settings,
where it was shown in \cite[Theorem 5]{A9} that the following theorem holds.
\begin{theorem}\label{propoo}
If $R$ is a Riemannian Jacobi-proportional algebraic curvature tensor which is not Clifford then $R$ is 2-root with multiplicities $8$ and $7$,
or it is 3-root with multiplicities $7$, $7$, and $1$.
\end{theorem}

At the end, we give generalizations to an indefinite case for two theorems from \cite[Theorems 3 and 4]{A9}.
As the proofs are based on the very same arguments as in \cite{A9}, we omit them.

\begin{theorem}\label{propo1}
Any two-root Jacobi-diagonalizable Osserman algebraic curvature tensor is Jacobi-proportional.
\end{theorem}

\begin{theorem}\label{propo2}
Any semi-Clifford algebraic curvature tensor is Jacobi-proportional.
\end{theorem}

\section*{Acknowledgements}

The authors are partially supported by the Ministry of Education, Science and Technological
Developments of the Republic of Serbia: grant number 451-03-68/2022-14/200104

%\section*{References}

%%\bibliography{mybibfile}

\end{document}